\documentclass[12pt, a4paper]{article}
\usepackage[T1]{fontenc}
\usepackage{amsmath,amsthm}
\usepackage{amssymb}
\usepackage{amscd}
\usepackage[latin1]{inputenc}
\usepackage{amsfonts}
\usepackage[matrix,arrow,curve]{xy}
\usepackage{graphicx}

\theoremstyle{plain}
\newtheorem{theo}{Theorem}
\newtheorem{lem}[theo]{Lemma}

\theoremstyle{definition}

\newtheorem{theosd}[theo]{Theorem}
\newtheorem{rem}[theo]{Remark}

\title{\Large A bound to kill the ramification over function fields
}
\author{\Large Alena Pirutka}

\begin{document}

\maketitle
\begin{abstract}
Let $k$ be a field of characteristic zero, let $X$ be a geometrically  integral $k$-variety of dimension $n$ and let $K$ be its field of fractions. Under the assumption that $K$ contains all $r^{\mathrm{th}}$ roots of unity for an integer $r$, we prove that, given an element $\alpha\in H^m(K,\mathbb Z/r)$, there exist $n^2$ functions $\{f_i\}_{,i=1,\ldots,n^2}$ such that $\alpha$ becomes unramified in $L=K(f_1^{1/r},\ldots,f_{n^2}^{1/r})$. 
\end{abstract}

\paragraph{1. Introduction.} Let $K$ be a field and let $\alpha\in \mathrm{Br}\,K$ be an element of order $r$. In \cite{S1}, \cite{S2}, Saltman proved that if $K$ is the function field of a $p$-adic curve and $(r,p)=1$, then $\alpha$ becomes trivial over an extension of $K$ of degree $r^2$. As a motivation for the question we consider in this paper, let us give a brief sketch of his arguments. Let us assume that $r$ is prime and that $K$ contains all $r^{\mathrm{th}}$ roots of unity. In fact, one can see that this case implies the general case.  We view $K$ as a function field of a regular, integral two-dimensional scheme $X$, projective over the spectrum of the ring of integers of a $p$-adic field.  Saltman then proved that one can find two functions $f_1, f_2\in K$ such that  $\alpha$ becomes unramified in $L=K(f_1^{1/r}, f_2^{1/r})$ with respect to any rank one discrete valuation ring  centered on $X$. This is sufficient to conclude, using the classical result that the Brauer group of a regular flat proper (relative) curve over the ring of integers of a $p$-adic field is trivial (cf. \cite{L}, \cite{T}).

Let us consider the case of higher dimensions, that is, assume that  $K$ is the field of fractions of an $n$-dimensional variety $X$, defined over a field $k$. Following Saltman's work, given a class $\alpha\in \mathrm{Br}\,K$, one may wonder if there is a bound $N$ depending only on $K$, such that we can kill the ramification of $\alpha$ with $N$ functions. Our main result (cf. theorem \ref{st}) gives an affirmative answer $N=n^2$ for $\alpha$ of order $r$ under the assumption that $K$ contains all $r^{\mathrm{th}}$ roots of unity. Our method also works for elements of  $H^m(K,\mathbb Z/r)$ and not only for $m=2$.
 
 
\paragraph*{Acknowledgements.} This work is a continuation of a discussion during the AIM workshop <<Deformation theory, patching, quadratic forms, and the Brauer group>> (Palo Alto, January 17-21, 2011). The author would like to thank the American Institute of Mathematics, Daniel Krashen and Max Lieblich for the organisation of this workshop and for their generous support.

\paragraph{2. Statement of the main result.}
Let $k$ be a field of characteristic zero. For $L$ a function field over $k$ containing all $r^{\mathrm{th}}$ roots of unity  we fix an isomorphism $\mu_r\stackrel{\sim}{\to}\mathbb Z/r$ of $Gal(\bar K/K)$-modules and we write
 \vspace{-0.1cm}
$$H_{nr}^m(L/k,\mathbb Z/r)=\bigcap_A\mathrm{Ker}[H^m(L, \mathbb Z/r)\stackrel{\partial_{A}}{\to}H^{m-1}(k_A, \mathbb Z/r)],$$
\vspace{-0.4cm}

\noindent where $A$ runs through all discrete valuation rings of rank one with $k\subset A$ and fraction field $L$. We denote by $k_A$ the residue field of $A$ and by $\partial_{A}$ the residue map.

\theosd\label{st}{\textit{Let $k$ be a field of characteristic zero. Let $X$ be an  integral $k$-variety of dimension $n$ and let $K$ be its field of fractions.  Let $r$ be an integer and assume that $K$ contains all $r^{\mathrm{th}}$ roots of unity. Let $\alpha$ be an element of $H^m(K,\mathbb Z/r)$. There exist $n^2$ functions $\{f_i\}_{,i=1,\ldots,n^2}$ such that $\alpha$ becomes unramified over $L=K(f_1^{1/r},\ldots,f_{n^2}^{1/r})$, that is,  we have  $\alpha_L\in H^m_{nr}(L/k, \mathbb Z/r)$. \\}}

We first prove two  lemmas.

\paragraph{3. Local description.}

In the case of dimension two, the following statement is due to Saltman (cf. \cite{S1} 1.2).
\lem\label{rl}{Let $k$ be an infinite field. Let $A$ be a local ring of a smooth $k$-variety and let $K$ be its field of fractions. Let $r$ be an integer prime to characteristic of $k$. Assume that $K$ contains all  $r^{\mathrm{th}}$ roots of unity and fix an isomorphism $\mu_r\stackrel{\sim}{\to}\mathbb Z/r$ of $Gal(\bar K/K)$-modules. Let $\alpha$ be an element of $H^m(K,\mathbb Z/r)$ ramified only at $s_1,\ldots, s_h$ forming a regular subsystem of parameters of the maximal ideal of $A$. Then $$\alpha=\alpha_0+\sum_{\emptyset\neq I\subset\{1,\ldots h\}} \alpha_I\cup s_I,$$ with $\alpha_0\in H^m(A,\mathbb Z/r)$, $\alpha_I\in H^{m-|I|}(A,\mathbb Z/r)$, and $s_I=\cup_{i\in I} (s_i)$, where we denote by $(s_i)$ the class of $s_i$ in $H^1(K,\mathbb Z/r)\simeq K^*/K^{*r}$.}
\proof{We proceed by induction on $h$ and $m$. Assume first $h=1$.  For $A$  a local ring of a smooth $k$-variety, with field of fractions $K$ and for $Y=Spec\,A$, we have an exact sequence due to Bloch and Ogus (cf. \cite{CTHK} 2.2.2)
\small
\begin{equation}\label{rec}
0\to H^m(A,\mathbb Z/r)\to H^m(K,\mathbb Z/r)\to\coprod_{x\in Y^{(1)}} H^{m-1}(\kappa(x), \mathbb Z/r)\to \coprod_{x\in Y^{(2)}} H^{m-2}(\kappa(x), \mathbb Z/r)\to\ldots\end{equation}
\normalsize
where the maps are induced by the residues.
Denote by $K(A/s_1)$ the field of fractions of $A/s_1$. As $\alpha$ is ramified only at $s_1$, we see from the sequence (\ref{rec})  that $\partial_{s_1}(\alpha)\in  H^{m-1}(K(A/s_1), \mathbb Z/r)$ is unramified. Hence,  from the sequence (\ref{rec})  for $A/s_1$, it comes from an element of $H^{m-1}(A/s_1, \mathbb Z/r)$. From Levine's conjecture (generalizing Bloch-Kato's conjecture proved by Rost and Voevodsky), proved by Kerz \cite{K} 1.2, any element of  $H^{m-1}(A/s_1, \mathbb Z/r)$ is a sum of cup products of units in $A/s_1$. In particular, any element of $H^{m-1}(A/s_1, \mathbb Z/r)$ lifts to $A$ : there exists an element $\alpha_1\in H^{m-1}(A, \mathbb Z/r)$ such that $\bar\alpha_1=\partial_{s_1}(\alpha)$. Hence $\alpha-\alpha_1\cup (s_1)$ is unramified, so it comes from $\alpha_0\in H^m(A,\mathbb Z/r)$, by (\ref{rec}) again.

If $m=1$, we have $\alpha=(s)$ for $s$ a function in $K$ and the result follows from the decomposition $s=u\prod_i s_i^{t_i}$ with $t_i\in \mathbb Z$ and $u\in A^*$.

Next, we assume the assertion for $(m-1,h-1)$ and $(m, h-1)$ and we prove it for $(m,h)$. From the sequence (\ref{rec}), $\partial_{s_1}(\alpha)\in  H^{m-1}(K(A/s_1), \mathbb Z/r)$ is ramified only at $\bar s_2,\ldots,\bar s_h$ where we denote by $\bar s_i$ the image of $s_i$ in $A/s_1$.  By induction, $\partial_{s_1}(\alpha)=\bar \alpha_1+\sum_{\emptyset\neq I\subset\{2,\ldots h\}} \bar \alpha_I\cup \bar s_I,$ where $\bar \alpha_1\in H^{m-1}(A/s_1,\mathbb Z/r)$, $\bar \alpha_I\in H^{m-1-|I|}(A/s_1,\mathbb Z/r)$, and $\bar s_I=\cup_{i\in I} (\bar s_i)$. As before, we deduce from \cite{K} 1.2 that all the $\bar \alpha_I$ and $\bar \alpha_1$ are sums of cup products of units in $A/s_1$ and so we can lift them to $ \alpha_I$ (resp. to $ \alpha_1$) on $A$. Now the element $\alpha-(\alpha_1+\sum_{\emptyset\neq I\subset\{2,\ldots h\}}  \alpha_I\cup  s_I)\cup (s_1)$ is ramified only at $s_2,\ldots,s_h$ and the lemma follows by induction.

\qed\\}

\paragraph{4. Divisor decomposition.}
\lem\label{ndiv}{Let $k$ be a field of characteristic zero and let $X$ be a  smooth projective  $k$-variety of dimension $n$. Let $D$ be a  divisor on $X$. There exists a  sequence of blowing-ups $f:X'\to X$ such that the support of   the total transform $f^*D$ is a simple normal crossing divisor which can be expressed a union of $n$ regular (but not necessarily connected) divisors of $X'$.}

\proof{By Hironaka, we may assume that $Supp(D)$ is a simple normal crossing divisor, which means that any irreducible component of $Supp(D)$ is smooth and that the fiber product over $X$ of any $c$ components of $Supp(D)$ is smooth and of codimension $c$. Let $G=(V,E)$ be the dual graph of $D$: \begin{itemize}
\item the vertices of $V$ correspond to  irreducible components $D_1,\ldots D_N$ of $D$
\item  the edge $(D_i, D_j)$ is in $E$ if the intersection $D_i\cap D_j$ is nonempty.
\end{itemize}

We say that we blow-up the edge $(D_i, D_j)$ if we change $X$ by the blow-up of the intersection $D_i\cap D_j$ (with reduced structure) and we change $G$ by the dual graph of the total transform of $D$, i.e. we add a vertex and corresponding edges. We write again $G=(V,E)$ for the modified graph.

We will show that after a finite sequence of blowing-ups $f:X'\to X$ of some edges  we may color the vertices of $G$ in $n$ colors so that for any edge $AB\in E$ the vertices $A$ and $B$ are of different colors.  Then $Supp(f^*D)$ is a simple normal crossing divisor and we  have  $Supp(f^*D)=\bigcup\limits_{i=1}^n F_i$ where $F_i$   is the (disjoint) union of components of $f^*D$ such that the corresponding vertex is of the $i^{\mathrm{th}}$ color. Hence $F_i$ are regular and the lemma follows.

If $n=2$  we may assume, after blowing-ups of some edges, that any cycle in $G$ has  even number of edges, which is sufficient to conclude.

Let us now assume that $n\geq 3$. We proceed by induction on the number $N$ of irreducible components of $D$. If $N\leq n$ the statement is clear. Assume it holds for $N$. Let $D$ be a divisor with $N+1$ components. By the induction hypothesis, after  blowing-ups of some edges,   we may assume that we may color all but the vertex $D_{N+1}$ of $G$ in $n$ colors as desired. We have $Supp(D)=\bigcup\limits_{i=1}^n F_i\;\cup D_{N+1}$ where $F_i$ is the union of components of $D$  of the $i^{\mathrm{th}}$ color. If $D_{N+1}$ doesn't intersect $F_i$ for some $i$ we color $D_{N+1}$ in $i^{\mathrm{th}}$ color. Hence we may assume that all the intersections $D_{N+1}\cap F_i$ are nonempty. By the same reason, we may assume that the intersection $F_2\cap F_3$ is nonempty. On the other hand, note that the intersection $\bigcap_i F_i\cap D_{N+1}$ is empty as $Supp(D)$ is a simple normal crossing divisor. We proceed by the following algorithm:

\begin{enumerate}
\item We first blow up all the edges $D_{j}D_{N+1}$ for all the components $D_j$ of $F_1$. Let us denote $E_1$ the union of all the exceptional divisors. This union is disjoint as the components of $F_1$ do not intersect. Note that  $E_1\cap F_2\cap\ldots \cap F_{n}=\emptyset$. Otherwise, we get a point in the intersection  $\bigcap_i F_i\cap D_{N+1}$ by projection. Moreover, there are no more edges between (the components of) $F_1$ and $D_{N+1}$ as the strict transforms of the corresponding divisors do not intersect.

\item Next, we blow up all the edges between $F_2$ and $F_3$ and we call $E_2$ the (disjoint) union of all new exceptional divisors.  Again, we have no more edges between $F_2$ and $F_3$ and also $E_2\cap E_1\cap F_4\cap\ldots \cap F_n=\emptyset$ (or $E_2\cap E_1$ is empty if $n=3$).
\item If $n=3$ we have the following picture:
\begin{center}
\includegraphics[scale=0.2]{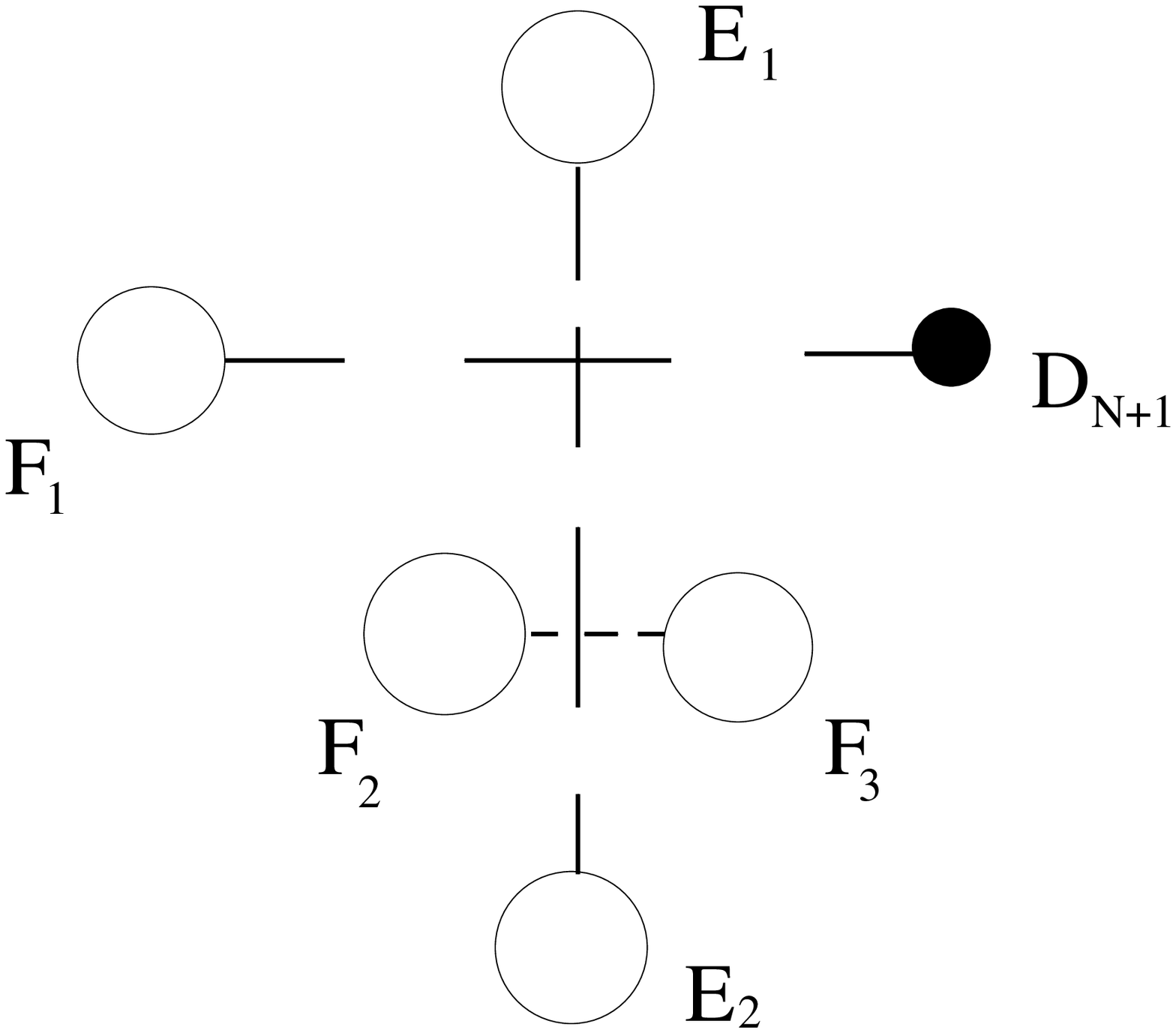}
\end{center}

 Here and in what follows the punctured line (for example, $F_2F_3$) means that there are no edges between components of corresponding groups (e.g. no edges between elements of $F_2\cup F_3$).

We color (all the vertices from) $F_1$ and $D_{N+1}$ in red, $E_1$ and $E_2$ in green and $F_2$ and $F_3$ in blue and this terminates the algorithm.
\item Assume that $n\geq 4$. We proceed until we get the group of exceptional divisors $E_{n-1}$ and then we go to step $6$. Suppose $3\leq i\leq n-1$ and we constructed $E_{i-2}$ and $E_{i-1}$ but no $E_i$. Suppose there are some edges between $E_{i-2}$ and $F_{i+1}$, otherwise we go to step $5$. We blow up all these edges and we call $E_i$ the (disjoint) union of all new exceptional divisors. We get no more edges between $E_{i-2}$ and $F_{i+1}$ and also $E_{i}\cap E_{i-1}\cap F_{i+2}\cap\ldots \cap F_n=\emptyset$.\\
\item If there are no edges between $E_{i-2}$ and $F_{i+1}$, we have the following picture:

\includegraphics[scale=0.5]{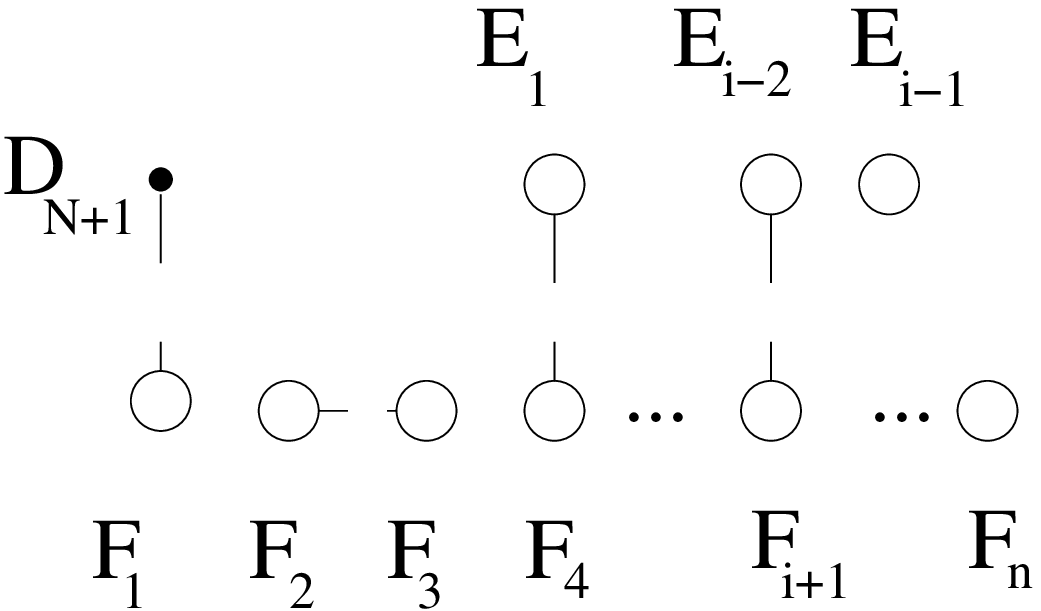}

    We color $F_1$ and $D_{N+1}$ in the first color, $F_2$ and $F_3$ in the second color, $E_1$ and $F_4$ in the third, ... , $E_{i-2}$ and $F_{i+1}$ in the $i^{\mathrm{th}}$-color, $E_{i-1}$ in  color $i+1$, and, finally, $F_{i+2}$, \ldots $F_n$ in colors $i+2, \ldots n$ respectively.

\item At this step, we have the following picture:

\includegraphics[scale=0.5]{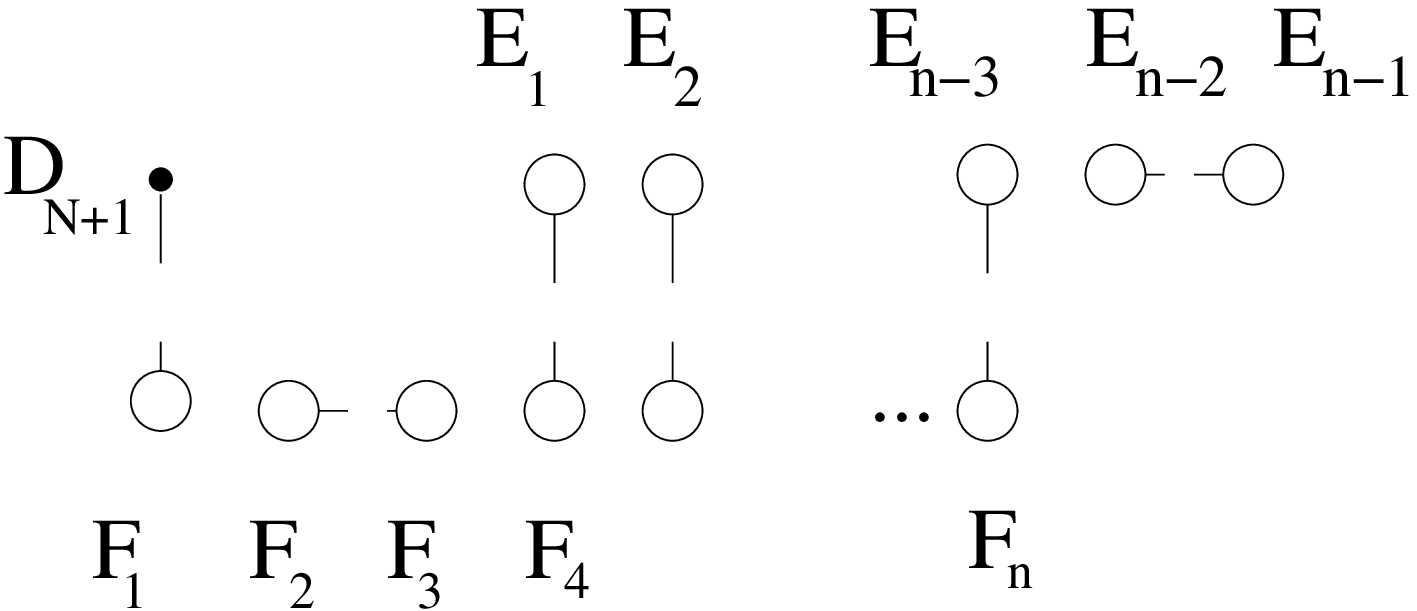}

    Moreover,  $E_{n-1}\cap E_{n-2}=\emptyset$ by construction. We color $F_1$ and $D_{N+1}$ in the first color, $F_2$ and $F_3$ in the second color, $E_1$ and $F_4$ in the third, ... , $E_{n-3}$ and $F_{n}$ in  color $n-1$,  $E_{n-1}$ and $E_{n-2}$ in  color $n$. This terminates the algorithm.\qed
\end{enumerate}
}

\paragraph{5. Proof of theorem \ref{st}.} By resolution of singularities, we may assume that $X$ is smooth.  By lemma \ref{ndiv}, we may assume that the ramification divisor $D=ram(\alpha)$ is a simple normal crossing divisor whose support is a union of $n$ regular divisors: $Supp\,D=\bigcup_{i=1}^n D_i$.   

For two divisors $G$ and $G'$ on $X$, with $G=\sum_{i=1}^{q} G_i$ where the  $G_i$ are irreducible divisors, we say that $G'$ is in \textit{general position} with $G$ if the support of $G'$  contains no generic point of any intersection $\bigcap_{i\in I} G_i$ for $I\subset\{1,\ldots,q\}$.

By a semilocal argument, we successively choose functions $f_1^j\in K$, $j=1,\ldots, n$, then $f_2^j\in K$, $j=1,\ldots, n$, \ldots, and then $f_n^j\in K$, $j=1,\ldots, n$, such that
$$div_X(f_i^j)=D_i+E_i^j$$
where $E_i^j$ are in general position with $D\cup\bigcup_{j'<j} Supp(E_i^{j'})$.

We claim that with this choice of $n^2$ functions $\alpha_L$ is unramified. Let $v$ be a discrete valuation on $L$ and let $x\in X$ be the point where the discrete valuation ring $R$ of $v$ is centered. We may assume that $x\in Supp\,D$, otherwise $\alpha$ is already unramified at $v$. From the construction, for any $i$, $D\cap\bigcap_{j=1}^nE_i^j=\emptyset$. Hence for any $1\leq i\leq n$ we can find $j_i$ such that $x\notin E_i^{j_i}$, which means that the corresponding local parameter $s_i$ of $D_i$ at $x$ is an $r^{\mathrm{th}}$ power in $K((f_i^j)^{1/r})$. Now the theorem follows from  lemma \ref{rl}, as any $s_I$ from the lemma is a cup product of $r^{\mathrm{th}}$ powers on $L$.
\qed\\

\rem{The bound $n^2$ is not sharp. For example, for $n=3$ one can kill all the ramification with four functions. Let us write $ram(\alpha)=D_1\cup D_2\cup D_3$ as in lemma \ref{ndiv}.  As in the proof of the theorem above, we take $f_i\in K$, $i=1,\ldots 4$, such that
\begin{gather*}
div(f_1)=D_1+D_2+D_3+E_1;\\
div(f_2)=D_1+D_2+E_2;\\
div(f_3)=D_2+D_3+E_3;\\
div(f_4)=D_1+2D_2+D_3+E_4.
\end{gather*}
and each $E_i$ is in general position with $ram(\alpha)\cup\bigcup_{i'<i} Supp(E_{i'})$.
Let $x$ be a center of a valuation $v$ on $L=K(f_i^{1/r})_{i=1,\ldots,4}$. We may assume that $x\in ram(\alpha)$. It is sufficient to see that if $x\in D_i$ then a local parameter of $D_i$ at $x$ can be expressed as a product of powers of the functions $f_i$.
\begin{enumerate}\item If $x\in X^{(1)}$ then $x$ lies on only one component $D_i$, which is thus defined by $f_1$.
 \item If $x\in X^{(2)}$ and if $x$ lies on two components $D_i$ and $D_j$, then $\frac{f_1}{f_3}$ defines $D_1$, $\frac{f_2f_3}{f_1}$ defines $D_2$, $\frac{f_1}{f_2}$ defines $D_3$. If $x$ lies on only one component $D_i$, then, by construction, $x\notin E_{i_1}\cup E_{i_2}$  for at least two indexes $1\leq i_1<i_2\leq 3$. By construction, $D_i$ is then defined at $x$ by at least one among the functions $f_{i_1}$ and $f_{i_2}$.
 \item Suppose that $x$ is a closed point of $X$. If $x\in D_1\cap D_2\cap D_3$, we use the same formulas as in the previous case. Next, suppose that $x$ lies on only two components of $ram(\alpha)$. Consider the case $x\in D_1\cap D_2$, the other cases are similar. By construction, $x$ lies on at most one component among $E_1, E_2,E_3$. Hence we see that if $x\notin E_2\cup E_3$ (resp. $x\notin E_1\cup E_3$, resp. $x\notin E_1\cup E_2$) then $D_1$ is defined by $\frac{f_2}{f_3}$ and $D_2$ is defined by $f_3$ (resp. by $\frac{f_1}{f_3}$ and by $f_3$, resp. by $\frac{f_1^2}{f_4}$ and by $\frac{f_4}{f_1}$).

     The last case is when $x$ lies on only one component of $ram(\alpha)$. Consider the case $x\in D_1$, the other cases are similar.  Then $x\notin E_1\cap E_2\cap E_4$ by construction.  Then $D_1$ is defined by $f_1$ (resp. by $f_2, f_4$) if $x$ does not lie on $E_1$ (resp. on $E_2, E_4$).
\end{enumerate}
}

\rem{By the same arguments as in the previous remark, if $r$ is prime to $2$ and $3$ and if $n=3$,  one can kill all the ramification with three functions $f_1,f_2,f_3$, such that
\begin{gather*}
div(f_1)=D_1+3D_2+3D_3+E_1;\\
div(f_2)=D_1+2D_2+D_3+E_2;\\
div(f_3)=D_1+D_2+2D_3+E_3.
\end{gather*}
and each $E_i$ is in general position with $ram(\alpha)\cup\bigcup_{i'<i} Supp(E_{i'})$.
}

\end{document}